# Optimality conditions and complete description of polytopes in combinatorial optimization

Alexey Antonov[*]

*University of Cambridge, UK. E-mail: aa2060@cam.ac.uk or aav.antonov@gmail.com

## Abstract

A combinatorial optimization problem (COP) has a finite groundset $E$ ($|E| = N$), a weight vector $c = (c^e: e \in E)$ and a family $T \in E$ of feasible subsets with objective to find $t \in T$ with maximal weight: $max\{\sum_{e \in t} c^e: t \in T\}$. Polyhedral combinatorics reformulates combinatorial optimization as linear program: $T$ is mapped into the set $X \in R^N$ of 0/1 incidence vectors and $c \in R^N$ is maximized over convex hull of $X$: $max\{cx: x \in conv(X)\}$. In theory, complementary slackness conditions for the induced linear program provide optimality conditions for the COP. However, for a general case, optimality conditions for combinatorial optimization have not been formulated analytically as for many problems complete description of the induced polytopes is available only as a convex hull of extreme points rather than a system of linear inequalities. Here, we formulate optimality conditions for a COP in general case: $x_k \in X$ is optimal if and only if $c \in cone(H_k)$ where $H_k = \{h \in V: hx_k \geq hx$ for any $x \in X\}$ and $V$ is a set of all -1/0/1 valued vectors in $R^N$. This provides basis to get, in theory, a complete description of a combinatorial polytope induced by any COP: all facet inducing inequalities for $conv(X)$ can be written as $hx \leq l$ where $h \in V$ and $l$ is integer. A vector $h \in V$ induces a nonredundent facet if and only if $h \in H_k^o \in H_k$ for at least one $x_k \in X$ (where $H_k^o = \{h \in H_k : h \notin cone(H_k \setminus \{h\})\}$) and $l = x_k h$.

## 1. Introduction and Notation

We consider a general linear combinatorial optimization problem for which we have a finite groundset $E$ ($|E| = N$), a weight vector $c = (c^e: e \in E)$ associated with each $e \in E$ and a family $T$ of feasible subsets of $E$. The objective is to find $t_* \in T$ that maximise the sum of weights associated with $t_*$. Thus we end up with the following Combinatorial Optimization Problem (Pulleyblank 1989):

$$max\{\sum_{e \in t} c^e: t \in T\} \quad \text{(COP)}$$

The COP is reformulated as linear programming problem by mapping each $t \in T$ into $R^N$ as (0/1) incidence vector $x$: $x^e = 1$ if $e \in t$ and $x^e = 0$ if $e \notin t$. We use $X$ to denote the set of all incidence vector $x \in R^N$ induced by $T$ and use $P^X = conv(X)$ to denote the corresponding combinatorial polytope (convex hull of $X$). Thus we can rewrite COP as Combinatorial Optimization Linear Program (Pulleyblank 1989, Aardal and Hoesel 1996):

$$max\{cx: x \in P^X\} \quad \text{(COLP)}$$

In the case of unique optimal solution COLP and COP are equivalent: the solution of COLP would be an extreme point of $P^X$ which is incidence vector $x_* \in X$ equivalent to the optimal feasible subset $t_* \in T$ for COP (Pulleyblank 1989).

Although polyhedral combinatorics reformulates combinatorial optimization as linear programs, optimality conditions for a general COP in analytical form (like complementary slackness conditions) have not been formulated. Analytical formulation of optimality conditions would require complete description of the induced combinatorial polytope as a system of linear inequalities (rather than a convex hull of incidence vectors) and is not known for many combinatorial problems. Here, for the first time (up to our knowledge), we provide analytical formulation (see theorem 2.1, 2.1b and 2.1c) of optimality conditions for a general linear COP. This, apparently, provides the basis for analytical formulation of a complete description (as a system of linear inequalities) of any combinatorial polytope (see theorem 3).

Notation that are used throughout the paper is summarised in table 1. In the table $y$, $b$, $z$ are points in $R^N$: $y \in R^N, b \in R^N, z \in R^N$ and $Z$ is a set of points in $R^N$: $Z \in R^N$. For example, we would frequently need to compare a product $yb$ with a product $zb$ for all points $z \in Z$ so we would regularly use notation 4 from table 1. We would also frequently check whether a point is inside or outside convex (conical) hull. Please note that further by decomposition (notation 5 and 6) we mean only nonnegative one ($y = \sum_{z_j \in Z} \gamma_j z_j$, $\gamma_j \geq 0$).

**Table 1.** Notation used throughout the paper

| ID | Notation | Meaning |
|---|---|---|
| 1 | yz | $\sum_{i=1}^{i=N} y^i z^i$ |
| 2 | $conv(Z)$ | convex hull of $Z$ |
| 3 | $cone(Z)$ | conical hull of $Z$ |
| 4 | $yb \geq Zb$ | $yb \geq zb$ for any $z \in Z$ |
| 5 | $y$ is decomposed by $Z$ | $\exists \gamma_j \geq 0$ so that $y = \sum_{z_j \in Z} \gamma_j z_j$ |
| 6 | decomposition of $y$ by $Z$ requires $z_k \in Z$ | if $y$ is decomposed by $Z$ and is not decomposed by $Z \setminus \{z_k\}$ |
| 7 | $H \in R^N$ | a hypercube with binary (0/1) vertexes |
| 8 | $H_V$ | vertexes of a unit hypercube (all 0/1 valued vectors in $R^N$) |
| 9 | $H^{\pm 1} \in R^N$ | convex hull (hypercube) of all -1/1 valued vectors in $R^N$. |
| 10 | $H(C)$ | a hypercube with a unit volume and -1/0/1 valued vertexes : $C \in E$ is indexing -1/0 valued vertexes of $H(C)$ $E/C$ is indexing 0/1 valued vertexes of $H(C)$ |

## 2. Optimality conditions in combinatorial optimization

To derive optimality conditions for COP we consider the set of all instances (SAI) of a given problem with a fixed groundset $E$ and a fixed family $T$ of feasible subsets. When both $E$ and $T$ are fixed SAI is defined by possible values of weight vector $c$. Due to the property of linear problems we can scale $c$ with any positive scalar (the order of solutions induced by $c$ is not affected by scaling of $c$) so that we can always ensure that $|c^e| < 1$ for each $e \in E$. Thus, without loss of generality SAI is reduced to a consideration of vector $c$ values from $N$-dimensional hypercube $H^{\pm 1} \in R^N$ which is the convex hull of all -1/1 valued vectors in $R^N$.

## 2.1 Weight vector: nonnegative case

First, for simplicity reasons, we would assume non-negativity of the weight vector: $c \geq 0$. In this case SAI is reduced to a unit hypercube $H \in R^N$: a hypercube with binary (0/1) vertexes. We would use $H$ to denote exactly a unit hypercube. We use $H_V$ to denote the set of all unit hypercube vertexes. Please note that $H_V$ is a set of all 0/1 valued vectors in $R^N$. In section 2.3 we would generalise our results by removing nonnegative constraints for the weight vector.

For each $x_k \in X$ we define a set $M_k = \{c \in H : x_k c \geq Xc\}$ of the weight vector values so that $x_k$ is optimal in COLP. From definition of $M_k$ it is obvious that if and only if $c \in M_k$ then $x_k$ is optimal. We will prove that $M_k = conv(H_k)$ where $H_k$ is the subset (defined by $x_k$ and $X$) of the unit hypercube vertexes only: $H_k = \{h \in H_V : x_k h \geq Xh\}$. Please note that by $x_k h \geq Xh$ we mean that $x_k h \geq xh$ for any $x \in X$ (see notation 4 from the table 1).

To illustrate our ideas consider a Maximum Spanning Tree problem with 3 nodes joined by 3 edges. In this trivial example (figure 1A) the set $X$ consist of 3 elements (possibilities to select 2 edges): $X = \{x_k = (0,1,1), x_j = (1,0,1), x_l = (1,1,0)\}$. The polytope $P^X$ is a convex hull of $X$ (blue triangle). Assuming $c \geq 0$ and properly scaled the set of all possible values for weight vector $c$ is a unit hypercube $H$ (in this case just a unit cube). Consider a point $x_k = (0,1,1)$ and a set $M_k = \{c \in H : x_k c \geq Xc\}$. It is easy to check that for $x_k = (0,1,1)$ to be optimal a weight vector $c$ must be in the polytope $M_k = conv(H_k)$ shown in figure 1B by yellow edges. The same polytopes $M_j$ and $M_l$ for points $x_j$ and $x_l$ are apparently convex hulls of corresponding sets $H_j$ and $H_l$ (not shown). The union of $M_k$, $M_j$ and $M_l$ is a unit cube $H$. It is also clear that if we consider the case $c \geq 0$ (without scaling) then for $x_k = (0,1,1)$ to be optimal we need just change convex with conical hull: $M_k = cone(H_k)$. Here we state that the same is related to any general linear COP (in case $c \in H$ the set of weight vector values so that an incidence vector $x \in X$ is optimal is defined as convex hull of a unit hypercube vertexes). We summarize this in theorem 2.1.

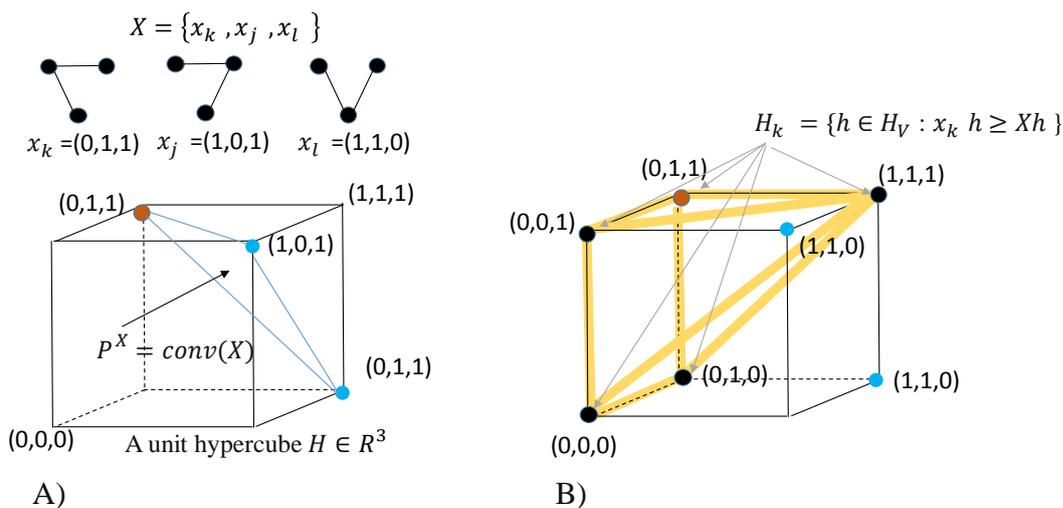

**Figure 1.** Maximum Spanning Tree (MST) problem with 3 nodes. **A)** If $c \geq 0$ then the set of all instances of MST problem is reduced to a cube of possible values of $c$ (assuming $c$ has been scaled properly). MST polytope $conv(X)$ is shown. **B)** The set $M_k = \{c \in H : x_k c \geq Xc\}$ for $x_k = (0,1,1)$ is shown by yellow edges as convex hull of cube vertexes: $M_k = conv(H_k)$. In case $c \in H$ it is easy to check that if and only if $c$ is inside yellow pyramid then $x_k$ is optimal.

**Theorem 2.1 (combinatorial optimization optimality conditions).** *Let $X \in R^N$ to be a set of incidence vectors induced by a COP with $c \geq 0$. Let $H \in R^N$ to be a unit hypercube and $H_V$ to be a set of H vertexes. Let $H_k \in H_V$ to be a subset of vertexes defined by $x_k \in X$ so that $H_k = \{h \in H_V : x_k\, h \geq Xh\}$. Then $x_k \in X$ is optimal if and only if $c \in cone(H_k)$:*

$$\exists\, \gamma_j \geq 0 \text{ so that } c = \sum_{h_j \in H_k} \gamma_j h_j \tag{1}$$

**Proof.** The proof of theorem 2.1 is directly follows from theorem 1 (see section 2.2). ∎

Please note that in theorem 2.1 we use $c \in cone(H_k)$ rather than $c \in conv(H_k)$ as we require in theorem $c \geq 0$ rather than $c \in H$. It is obvious that if $c \in conv(H_k)$ then $\alpha c \in cone(H_k)$ for any $\alpha \geq 0$.

## 2.2 Proving theorem 2.1

It is clear that the set $X$ of all incidence vectors for any COP is a subset of $H_V : X \in H_V$. Therefore, theorem 1 provides the basis for theorem 2.1:

**Theorem 1.** *Let H to be a unit hypercube and $H_V$ to be a set of H vertexes. Let $X \in H_V$ to be a subset of $H_V$. For any vertex $x_k \in X$ the set of points $M_k = \{x \in H : x_k\, x \geq Xx\}$ is a convex hull of $H_k = \{h \in H_V : x_k\, h \geq Xh\}$:*
$$M_k = conv(H_k).$$

To prove Theorem 1 we will introduce several sets defined by a unit hypercube vertexes ($x_k$, $x_k$ and $x_l$). Initially we introduce the set $B_{kj} = \{x \in H : x_k\, x = x_j\, x\}$ which is a boundary between any two vertexes $x_k$ and $x_j$. Lemma 1 states that such a boundary inside a unit hypercube is a convex hull of hypercube vertexes $H_{kj} = \{h \in H_V : x_k\, h = x_j\, h\}$ (see figure 2A). The set $B_{kjl} = \{x \in H : x_k\, x = x_j\, x = x_l\, x\}$ extends definition of $B_{kj}$ for 3 vertexes $x_k$, $x_j$ and $x_l$. Lemma 2 states that $B_{kjl} = conv(H_{kjl})$ where $H_{kjl} = \{h \in H_V : x_k\, h = x_j\, h = x_l\, h\}$. It is obvious that $B_{kjl} \in B_{kj}$ as well as $H_{kjl} \in H_{kj}$. The set $B_{kj}^l = \{x \in H : x_k\, x = x_j\, x \geq x_l\, x\}$ (lemma 3) is also defined by 3 vertexes $x_k$, $x_j$ and $x_l$ but represents a different subset of $B_{kj}$ (see figure 2B) and also is a convex hull of hypercube vertexes $H_{kj}^l = \{h \in H_V : x_k\, h = x_j\, h \geq x_l\, h\}$. Finally, a set $B_{kj}^Y = \{x \in H: x_k\, x = x_j\, x \geq Yx\}$ is an extension of $B_{kj}^l$ ( a set of vertexes $Y$ vs. a single $x_l$ ). Lemma 4 states that $B_{kj}^Y$ is a convex hull of hypercube vertexes $H_{kj}^Y = \{h \in H_V : x_k\, h = x_j\, h \geq Yh\}$. So we have the following relations between sets:

$$H_{kj}^Y \in H_{kj}^l \in H_{kj} \in H_k \quad \text{for any hypercube vertexes } x_k,\ x_j\ , x_l\ \in Y$$

$$B_{kj}^Y \in B_{kj}^l \in B_{kj} \in B_k \quad \text{for any hypercube vertexes } x_k,\ x_j\ , x_l\ \in Y$$

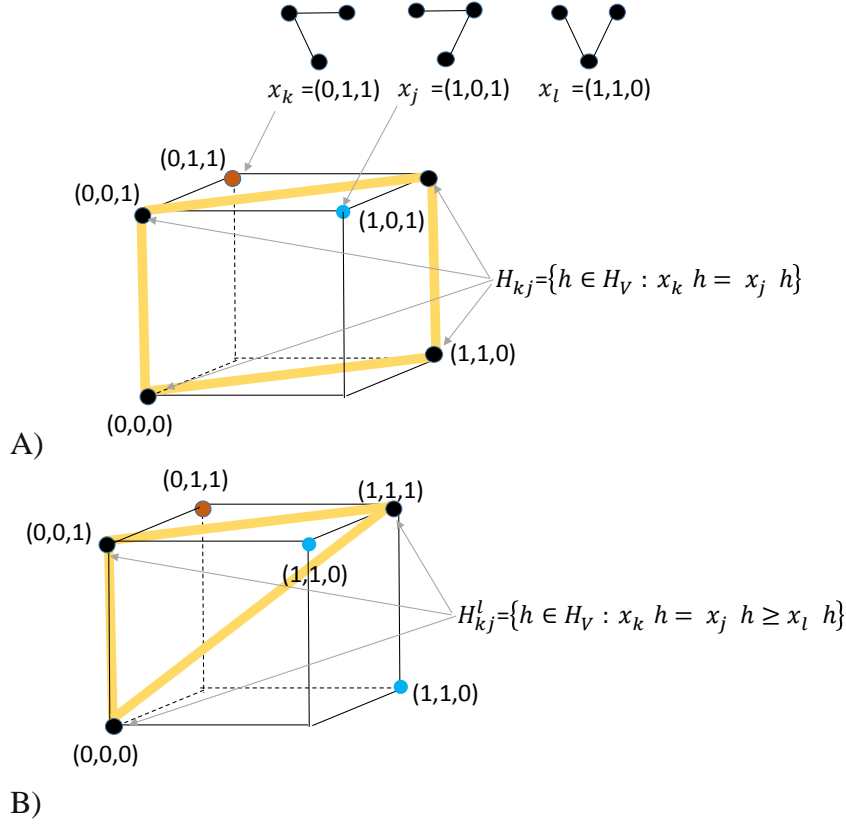

**Figure 2.** Sets ($B_{kj}$ and $B_{kj}^l$) used to prove theorem 1 for MST problem with 3 nodes. **A)** A set $B_{kj} = \{x \in H : x_k\ x = x_j\ x\}$ for $x_k = (0,1,1)$ and $x_j = (1,1,0)$ is shown by yellow edges: the set of points $B_{kj}$ is defined as convex hull of a unit hypercube vertexes $H_{kj} = \{h \in H_V : x_k\ h = x_j\ h\}$. **B)** A subset $B_{kj}^l = \{x \in H : x_k\ x = x_j\ x \geq x_l\ x\}$ is shown with yellow edges. $B_{kj}^l$ is defined by $H_{kj}^l = \{h \in H_V : x_k\ h = x_j\ h \geq x_l\ h\}$.

**Lemma 1.** *The subset* $B_{kj} = \{x \in H : x_k\ x = x_j\ x\}$ *for any two vertexes* $x_k \in H_V$ *and* $x_j \in H_V$ *is a convex hull of* $H_{kj} = \{h \in H_V : x_k\ h = x_j\ h\}$:
$$B_{kj} = conv(H_{kj}).$$

**Proof.** First it is obvious that if $x \in conv(H_{kj})$ then $x_k\ x = x_j\ x$ as $x$ is decomposed by $H_{kj}$: $x = \sum_{h_j \in H_{kj}} \gamma_j h_j$, $\gamma_j \geq 0$ and for each $h_j \in H_{kj}$ we have $x_k\ h_j = x_j\ h_j$. Next we prove that any point outside $conv(H_{kj})$ is outside $B_{kj}$.

The set $B_{kj}$ is an intersection of a hyperplane $P_{kj} = \{p \in R^N : (x_k - x_j)p = 0\}$ and a unit hypercube $H$: $B_{kj} = H \cap P_{kj}$. Any vertex $h \in H_{kj} \in B_{kj}$ is an extreme point of $B_{kj}$ (as $h$ is an extreme point of $H$ and $B_{kj} \in H$). The set $B_{kj}$ will have different (not $H_{kj}$) extreme points only if there is a 1-dimentional edge $s$ of the hypercube $H$ so that a hyperplane $P_{kj}$ intersects with $s$ at a single point $s_o \in B_{kj}$:
$$s \cap P_{kj} = s_o\ ,\ \ s_o \notin H_{kj}$$

In fact this is not possible due to the binary nature of the unit hypercube vertexes.

Let us recall that two vertexes $h_1$ and $h_0$ of a unit hypercube are connected by 1-dimentinal edge $s$ if they have only one different element. Let us denote $h_1$ and $h_0$ in a way that

$$(x_k - x_j)h_1 \geq (x_k - x_j)h_0 \qquad (11.1)$$

We can always do this because $h_1$ and $h_0$ is 0/1 valued and have only one different element while $(x_k - x_j)$ is -1/0/1 valued. Thus, for any two vertexes $h_1$ and $h_0$ that are connected by 1-dimentinal edge is true:

$$0 \leq (x_k - x_j)h_1 - (x_k - x_j)h_0 \leq 1 \qquad (11.2)$$

Now assume that edge $s$ intersects with $P_{kj}$ at a single point $s_o \notin H_{kj}$. This means that $h_1$ and $h_0$ are located at different halve spaces separated by hyperplane $P_{kj}$ and therefore (assuming we denote $h_1$ and $h_0$ in a way 11.1):

$$(x_k - x_j)h_1 > 0 \qquad (11.3)$$
$$(x_k - x_j)h_0 < 0 \qquad (11.4)$$

However, this is not possible. If (11.3) is true then $(x_k - x_j)h_1 \geq 1$ as $(x_k - x_j)h_1 > 0$ and integer. In similar way we get that $(x_k - x_j)h_0 \leq -1$. This means that $(x_k - x_j)h_1 - (x_k - x_j)h_0 \geq 2$ what contradicts to (11.2) and, therefore, (11.3, 11.4) is false if (11.2) is true. This means that hyperplane $P_{kj}$ can not intersects with 1-dimentional edges at a single point different from $H$ vertex and, thus, the set of extreme points of $B_{kj} = H \cap P_{kj}$ consist of only $H$ vertexes which, in turn, proves that $B_{kj} = conv(H_{kj})$. ∎

**Lemma 2.** *The subset $B_{kjl} = \{x \in H : x_k\,x = x_j\,x = x_l\,x\}$ for any $x_k \in H_V$, $x_j \in H_V$ and $x_l \in H_V$ is a convex hull of $H_{kjl} = \{h \in H_V : x_k\,h = x_j\,h = x_l\,h\}$:*
$$B_{kjl} = conv(H_{kjl})$$

**Proof.** It is trivial to prove that $conv(H_{kjl}) \in B_{kjl}$. Next we prove that any point outside $conv(H_{kjl})$ is outside $B_{kjl}$.

Assume there is point $z \in B_{kjl}$: $x_k\,z = x_j\,z = x_l\,z$ and $z \notin conv(H_{kjl})$. From lemma 1 it follows that $z \in conv(H_{kj})$. Thus we can decompose $z$ by $H_{kj}$:

$$z = \sum_{h_i \in H_{kj}} \gamma_i h_i \quad , \quad \gamma_i \geq 0$$

By definition $H_{kjl} \in H_{kj}$ so we can rewrite decomposition of $z$:

$$z = \sum_{h_i \in H_{kjl}} \gamma_i h_i + \sum_{h_i \in (H_{kj} \setminus H_{kjl})} \beta_i h_i \quad , \quad \gamma_i \geq 0 \text{ and } \beta_i \geq 0$$

For any $h_i \in H_{kjl}$ we have $x_k\, h_i = x_j\, h_i = x_l\, h_i$ while for any $h_i \in (H_{kj}\setminus H_{kjl})$ we have $x_k\, h_i = x_j\, h_i > x_l\, h_i$. This means that either each $\beta_i = 0$ or $x_k\, z = x_j\, z > x_l\, z$. In the first case ($\beta_i = 0$) we have $z = \sum_{h_i \in H_{kjl}} \gamma_i h_i$, $\gamma_i \geq 0$ which contradicts to $z \notin conv(H_{kjl})$ while the second case ($x_k\, z = x_j\, z > x_l\, z$) contradicts to $z \in B_{kjm}$. ∎

**Lemma 3.** The subset $B_{kj}^l = \{x \in H : x_k\, x = x_j\, x \geq x_l\, x\}$ for any $x_k \in H_V$, $x_j \in H_V$ and $x_l \in H_V$ is a convex hull of $H_{kj}^l = \{h \in H_V : x_k\, h = x_j\, h \geq x_l\, h\}$:
$$B_{kj}^l = conv(H_{kj}^l)$$

**Proof.** It is trivial to prove that $(H_{kj}^l) \in B_{kj}^l$. Next we prove that any point outside $conv(H_{kj}^l)$ is outside $B_{kj}^l$. For this we would need to use the relations between introduced subsets:

$$H_{kjl} \in H_{kj}^l \in H_{kj} \text{ and } B_{kjl} \in B_{kj}^l \in B_{kj} \tag{1.3.0}$$

These relations follows from the subsets definitions.

Assume that there is a point $z \in B_{kj}^l$:
$$x_k\, z = x_j\, z \geq x_l\, z \tag{1.3.1}$$

and

$$z \notin conv(H_{kj}^l) \tag{1.3.2}$$

If $x_k\, z = x_j\, z = x_l\, z$ then $z \in B_{kjl}$ and from lemma 2 $z \in conv(H_{kjl})$. From (1.3.0) we have $H_{kjl} \in H_{kj}^l$ and, hence, $z \in conv(H_{kj}^l)$ which is contradiction to 1.3.2. Therefore, we need to consider only the case $x_k\, z = x_j\, z > x_l\, z$.

From $x_k\, z = x_j\, z$ it follows that $z \in B_{kj}$ and from lemma 1 it follows that $z \in conv(H_{kj})$. Thus we can decompose $z$ using vertexes from $H_{kj}$:
$$z = \sum_{h_i \in H_{kj}} \gamma_i h_i, \quad \gamma_i \geq 0$$

As $H_{kj}^l \in H_{kj}$ we can rewrite decomposition as

$$z = \sum_{h_i \in H_{kj}^l} \gamma_i h_i + \sum_{h_i \in (H_{kj}\setminus H_{kj}^l)} \beta_i h_i \text{ where } \gamma_i \geq 0 \text{ and } \beta_i \geq 0 \tag{1.3.3}$$

Let us select $h_* \in (H_{kj}\setminus H_{kj}^l)$ so that $\beta_* > 0$ (if there is no such $\beta_*$ then $z = \sum_{h_i \in H_{kj}^l} \gamma_i h_i$, $\gamma_i \geq 0$ and, therefore, $z \in conv(H_{kj}^l)$). For $h_*$ (due to $h_* \in (H_{kj}\setminus H_{kj}^l)$) we have $x_k\, h_* = x_j\, h_* < x_l\, h_*$ and for $z$ (due to $z \in B_{kj}^l$) we have $x_k\, z = x_j\, z > x_l\, z$ which means that there is $\mu: 0 < \mu < 1$ so that for point $z^o = \mu h_* + (1-\mu)z$ is true:

$$x_k\, z^o = x_j\, z^o = x_l\, z^o \tag{1.3.4}$$

We can express $h_*$ as linear combination of $z^o$ and $z$:

$$h_* = (1/\mu)z^o - ((1-\mu)/\mu)z \tag{1.3.5}$$

From (1.3.4) it follows that $z^o \in B_{kjl}$ and, therefore, by lemma 2 we have $z^o \in conv(H_{kjl})$ and thus $z^o \in conv(H_{kj}^l)$ as $H_{kjl} \in H_{kj}^l$ which means that we can decompose $z^o = \sum_{h_i \in H_{kj}^l} a_i h_i$ with $a_i \geq 0$ and rewrite (1.3.5):

$$h_* = \left(\frac{1}{\mu}\right)\sum_{h_i \in H_{kj}^l} a_i h_i - \left(\frac{1-\mu}{\mu}\right) z, a_i \geq 0 \tag{1.3.6}$$

Substituting (1.3.6) into (1.3.3) we will eliminate $\beta_* h_*$ (i.e. $\beta_* = 0$):

$$z = \sum_{h_i \in H_{kj}^l} \gamma_i^{new} h_i + \sum_{h_i \in (H_{kj} \setminus H_{kj}^l)} \beta_i^{new} h_i, \quad \gamma_i^{new} \geq 0, \quad \beta_i^{new} \geq 0 \text{ and } \beta_*^{new} = 0$$

where the new coefficients are:

$$\gamma_i^{new} = (\mu/(\mu + \beta_*(1-\mu)))\gamma_i + (\beta_*/(\mu + \beta_*(1-\mu))) a_i$$

$$\beta_i^{new} = (\mu/(\mu + \beta_*(1-\mu)))\beta_i$$

and $a_i \geq 0$, $\gamma_i \geq 0$, $\beta_* > 0$, $1 > \mu > 0$.

It is easy to see that $\gamma_i^{new} \geq 0$. It is important that $\beta_i^{new} > 0$ only if $\beta_i > 0$. Thus l.3.3 can be rewritten without $h_*$ ($\beta_*^{new} = 0$ and without adding novel $h_i \in (H_{kj} \setminus H_{kj}^l)$ with $\beta_i > 0$). Repeating similar procedure for every $\beta_i > 0$ we will finally rewrite 1.3.3 so that for any $h_i \in (H_{kj} \setminus H_{kj}^l)$ the corresponding $\beta_i = 0$ which means that $z$ is decomposed by $H_{kj}^l$. Thus, if 1.3.1 is true then $\exists \gamma_i \geq 0$ so that $z = \sum_{h_i \in H_{kj}^l} \gamma_i h_i$ and 1.3.2 is false. ∎

Trivial consequence from lemma 3 is that if decomposition of $z$ using vertexes $H_{kj}$ requires a nonzero input from any vertex $h_o \notin H_{kj}^l$ then $z \notin B_{kj}^l$ as in this case $z \notin conv(H_{kj}^l)$.

**Lemma 4.** *The subset $B_{kj}^Y = \{x \in H : x_k\ x = x_j\ x \geq Yx\}$ for any $x_k \in H_V$ and $x_j \in H_V$ and a set $Y \in H_V$ is a convex hull of $H_{kj}^Y = \{h \in H_V : x_k\ h = x_j\ h \geq Yh\}$:*

$$B_{kj}^Y = conv(H_{kj}^Y)$$

**Proof.** It is trivial to prove that $conv(H_{kj}^Y) \in B_{kj}^Y$. Next we prove that any point outside $conv(H_{kj}^Y)$ is outside $B_{kj}^l \cap B_{kj}^t$.

Both sets $B_{kj}^Y$ and $H_{kj}^Y$ can be presented as intersections (it follows from the sets definitions):

$$B_{kj}^Y = \bigcap_{x_l \in Y} B_{kj}^l \text{ and } H_{kj}^Y = \bigcap_{x_l \in Y} H_{kj}^l$$

Assume that there is a point $z$ so that:

$$z \in B_{kj}^Y = \bigcap_{x_l \in Y} B_{kj}^l \tag{1.4.1}$$

and

$$z \notin conv(H_{kj}^Y) = conv(\bigcap_{x_l \in Y} H_{kj}^l) \tag{1.4.2}$$

For each $x_l \in Y$ we have $z \in B_{kj}^l$ and by lemma 3 we can decompose $z = \sum_{h_i \in H_{kj}^l} \gamma_i h_i$ with $\gamma_i \geq 0$. From (1.4.2) it follows that there is at least one pair $x_l \in Y$ and $x_t \in Y$ so that there is a vertex $h_o \in H_{kj}^l$ and $h_o \notin H_{kj}^t$ and $h_o$ is required for decomposition of $z$ with vertexes $H_{kj}$. Therefore, based on the consequence of lemma 3 we have $z \notin B_{kj}^t$ and therefore $z \notin \bigcap_{x_l \in Y} B_{kj}^l = B_{kj}^Y$ as $\bigcap_{x_l \in Y} B_{kj}^l \in B_{kj}^t$. ∎

**Proof of Theorem 1.** It is trivial to prove that $conv(H_k) \in M_k$. Next we prove that any point outside $conv(H_k)$ is outside $M_k$.

Assume that there is point $z \in M_k$ so that:

$$x_k \, z > Xz \tag{T1.1a}$$

and

$$z \notin conv(H_k) \tag{T1.1b}$$

Please note that we consider the case $x_k \, z > Xz$ rather than $x_k \, z \geq Xz$ because if there is $x_j \in X$ so that $x_k \, z = x_j \, z \geq Xz$ then $z \in B_{kj}^X$ and, therefore, by lemma 4, $z \in conv(H_{kj}^X)$. And due to $H_{kj}^X \in H_k$ it automatically implies that $z \in conv(H_k)$. Therefore, further we need to consider only the case $x_k \, z > Xz$.

As $z \in M_k \in H$ we can decompose $z$ by $H_V$:

$$z = \sum_{h_i \in H_k} \gamma_i h_i + \sum_{h_i \in (H_V \setminus H_k)} \beta_i h_i \quad \text{where } \gamma_i \geq 0 \text{ and } \beta_i \geq 0 \tag{T1.2}$$

Let us consider vertex $h_* \in (H_V \setminus H_k)$ so that $\beta_* > 0$ (if there is no such vertex $h_*$ then $z = \sum_{h_i \in H_k} \gamma_i h_i$ with $\gamma_i \geq 0$ and therefore $z \in conv(H_k)$ which makes T1.1b false). For vertex $h_*$ there is at least one $x_j \in X$ so that $x_k \, h_* < x_j \, h_*$ (otherwise $h_* \in H_k$). Let $X_l$ to be all such vertexes:

$$X_l = \{x_l \in X : x_k \, h_* < x_l \, h_*\}$$

On the other side we have

$$x_k \, z > x_l \, z \text{ for each } x_l \in X_l$$

which means that for each $x_l \in X_l$ there is a point $z^l$:

$$z^l = \mu_l \, h_* + (1 - \mu_l)z, \; 1 > \mu_l > 0$$

so that

$$x_k \, z^l = x_l \, z^l.$$

Let us select $x_o \in X_l$ with the smallest $\mu_o : \mu_o \leq \mu_l$ for each $x_l \in X_l$. We denote

$$z^o = \mu_o \, h_* + (1 - \mu_o)z \tag{T1.3}$$

For $z^o$ and $x_o$ we have

$$x_k \, z^o = x_o \, z^o \tag{T1.4}$$

It is important to notice that $x_k \, z^o \geq x_l \, z^o$ for each $x_l \in X_l$ (due to $\mu_o \leq \mu_l$). For all points $x_j \in X \setminus X_l$ we have $x_k \, z^o > x_j \, z^o$ as $x_k \, z > x_j \, z$ and $x_k \, h_* \geq x_j \, h_*$. Hence for $z^o$ we have $x_k \, z^o \geq xz^o$ for any $x \in X$:

$$x_k \, z^o \geq Xz^o \tag{T1.5}$$

From T1.4 and T1.5 it follows that $x_k \, z^o = x_o \, z^o \geq Xz^o$ which means that $z^o \in B_{ko}^X$ and from lemma 4 it follows that $z^o \in conv(H_{ko}^X)$. From $H_{ko}^X \in H_{ko} \in H_k$ it follows that

$$z^o \in conv(H_k) \tag{T1.5}$$

We can express $h_*$ from (T1.3):

$$h_* = (1/\mu_o)z^o - ((1 - \mu_o)/\mu_o)z \tag{T1.6}$$

Due to T1.5 we can decompose $z^o = \sum_{h_i \in H_k} a_i h_i$, $a_i \geq 0$ using only $H_k$ vertexes and therefore

$$h_* = (1/\mu_o) \sum_{h_i \in H_k} a_i h_i - ((1 - \mu_o)/\mu_o)z, \; a_i \geq 0 \tag{T1.7}$$

Substituting (T1.7) into (T1.2) we will eliminate $\beta_* h_*$ in decomposition of $z$:

$$z = \sum_{h_i \in H_k} \gamma_i^{new} h_i + \sum_{h_i \in (H_V \setminus H_k)} \beta_i^{new} h_i, \; \gamma_i^{new} \geq 0, \; \beta_i^{new} \geq 0 \text{ and } \beta_*^{new} = 0 \tag{T1.8}$$

where new coefficients are:

$$\gamma_i^{new} = (\mu_o/(\mu_o + \beta_*(1 - \mu_o)))\gamma_i + (\beta_*/(\mu_o + \beta_*(1 - \mu_o))) \, a_i$$

$$\beta_i^{new} = (\mu_o / (\mu_o + \beta_*(1 - \mu_o)))\beta_i$$

and $a_i \geq 0$, $\gamma_i \geq 0$, $\beta_* > 0$, $1 > \mu > 0$.

For all coefficients $\gamma_i^{new}$ we have $\gamma_i^{new} \geq 0$. It is important to note that $\beta_i^{new} > 0$ only if $\beta_i > 0$. Therefore, a new decomposition (T1.8) of $z$ does have a zero input ($\beta_*^{new} = 0$) from $h_*$ while no new positive inputs from $h_i \in (H_V \setminus H_k)$ are added. We can repeat the same procedure for each remaining $\beta_i > 0$ for vertexes $h_i \in (H_V \setminus H_k)$ in T1.8 and, thus, making all $\beta_i^{new}$ for vertexes $h_i \in (H_V \setminus H_k)$ to be zero. So this proves that if T1.1a is true then we can always decompose $z$ by $H_k$ vertexes only: $\exists\ \gamma_i \geq 0$ so that $z = \sum_{h_i \in H_k} \gamma_i h_i$. Therefore, if $z \in M_k$ then $z \in conv(H_k)$ and $M_k = conv(H_k)$ ∎

## 2.3 Weight vector: general case

Now we remove nonnegativity constraint from the weight vector and consider several possible cases. As already been mentioned, in general case, $c \in R^N$ and we can consider $c \in H^{\pm 1}$ (as for linear problems we can always rescale weight vector to ensure this). It is easy to notice that hypercube $H^{\pm 1} \in R^N$ consist of $2^N$ hypercubes with a unit volume including the unit hypercube $H$ (0/1 valued vertexes). We denote such hypercubes $H(C)$ (figure 3) where $C \in E$. For each vertex of $H(C)$ hypercube components indexed by $C$ are -1/0 valued while components indexed by $E \setminus C$ are 0/1 valued. Each vertex $h(C) \in H(C)$ is induced from vertex $h \in H$ by the following transformation:

$$h(C) = \begin{cases} h^i - 1, & i \in C \\ h^i, & i \in E \setminus C \end{cases} \qquad (2.3.1)$$

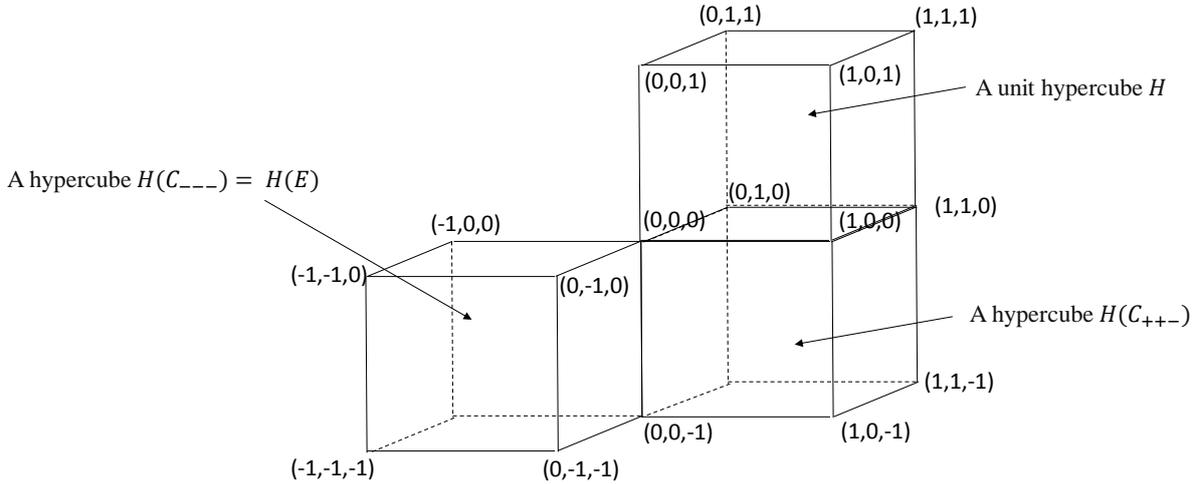

**Figure3.** A hypercube $H^{\pm 1} \in R^N$ consist of $2^N$ hypercubes $H(C)$ with a unit volume. Any $H(C)$ can be induced from a unit hypercube $H = H(C = empty)$ using 2.3.1 transformation for each vertex from $H$.

It is easy to see the link between $H(C)$ and the weight vector: if $c \in H(C)$ then $C$ is a set of indexes where the weight vector components are negative: $C = \{e \in E : c^e < 0\}$. First we consider case $c \in H(C)$ and state that nothing will change in theorem 2.1 except that the set $H_V$ would be substituted with a set $H_V(C)$ where $H_V(C)$ is a set of $H(C)$ vertexes. Finally, we consider a general case $c \in H^{\pm 1}$. We state that in this case one need to substitute $H_V$ in theorem 2.1 with a set $\bigcup_{C \in E} H_V(C)$ which is a set $V$ of all -1/0/1 valued vectors in $R^N$.

We have proved Theorem 1 for a unit hypercube $H$ with 0/1 valued vertexes. To prove it we use only two properties of $H$. First, throughout the proof we used the property that $H$ is convex and thus any $x \in H$ can be decomposed as nonnegative linear combination of $H$ vertexes. Second, in lemma 1 we used the property of a unit hypercube $H$ that each 1-dimensional edge of $H$ is formed by 2 vertexes that have only one different component and that the difference between integer values (in this case 0 and 1) is 1. But any other hypercube with a unit length of 1-dimensional edge and integer valued vertexes has the same properties. Thus we can extend theorem 1 by substituting a unit hypercube with any hypercube that have a unit volume and integer valued vertexes:

**Theorem 1b.** *Let $H^1$ to be a hypercube with a unit volume and integer valued vertexes. Let $H_V^1$ to be the set of all $H^1$ vertexes and $X \in H_V^1$. For any vertex $x_k \in X \in H_V^1$ the set of $H^1$ points $M_k = \{x \in H^1 : x_k x \geq Xx\}$ is a convex hull of $H_k = \{h \in H_V^1 : x_k h \geq Xh\}$:*
$$M_k = conv(H_k)$$

Let us remark that theorem 1b is valid for any hypercube $H(C)$ as each $H(C)$ has a unit volume and integer (-1/0/1) valued vertexes. Now let us denote $H_V(C)$ the set of $H(C)$ vertexes. For a subset $X \in H_V$ we denote $X(C)$ a subset of $H_V(C)$ induced by $X$ using transformation (2.3.1):

$$x_k(C) = \begin{cases} x_k^i - 1, & i \in C \\ x_k^i, & i \in E \setminus C \end{cases}, \text{ where } x_k \in X \in H_V \text{ and } x_k(C) \in X(C) \in H_V(C)$$

**Lemma 2.3.1** *For any $x_k \in H_V$ and $x_j \in H_V$ and for any point $x \in H(C)$ we have:*

$$x_k x - x_k(C)x = x_j x - x_j(C)x \tag{l2.3.1}$$

**Proof.** For any vertex $h \in H^V$ of a unit hypercube $H$ and vertex $h(C) \in H_V(C)$ of a unit hypercube $H(C)$ and any $x \in H(C)$ we have:

$$hx - h(C)x = \sum_{i \in (E \setminus C)}(h_k^i x^i - h_k^i x^i) + \sum_{i \in C}(h_k^i x^i - (h_k^i - 1)x^i) = \sum_{i \in C} x^i$$

Therefore, and the difference $hx - h(C)x = \sum_{i \in C} x^i$ does not depend on $h$. For each $x_k \in H_V$ and $x_j \in H_V$ we have $x_k x - x_k(C)x = \sum_{i \in C} x^i = x_j x - x_j(C)x$ which proves (l2.3.1). ∎

As a consequence of lemma 2.3.1 we have $x_k x - x_j x = x_k(C)x - x_j(C)x$. This means that any equality/inequality valid for $x_k \in H_V$, $x_j \in H_V$ and $x \in H(C)$ is valid for $x_k(C) \in H_V(C)$, $x_j(C) \in H_V(C)$ and $x \in H(C)$: if $x_k x \geq x_j x$ then the same is true for $x_k(C)x \geq x_j(C)x$ (and vice versa).

Consider a set $M_k(C) = \{x \in H(C) : x_k\, x \geq Xx\}$ defined for each $x_k \in X \in H_V$. Consider analogous set $M'_k(C) = \{x \in H(C) : x_k(C)x \geq X(C)x\}$ defined for each $x_k(C) \in X(C) \in H_V(C)$. Similar we define $H_k(C) = \{h \in H_V(C) : x_k\, h \geq Xh\}$ and $H'_k(C) = \{h \in H_V(C) : x_k(C)h \geq X(C)h\}$.

**Lemma 2.3.2** $M_k(C) = M'_k(C)$ and $H_k(C) = H'_k(C)$.
**Proof.** The proof of lemma 2.3.2 follows directly from lemma 2.3.1. ∎

From Theorem 1b and lemma 2.3.2 it follows that theorem 1 is valid if we substitute a unit hypercube $H$ with a hypercube $H(C)$ and modify definition of $M_k$ and $H_k$ sets to be subsets of $H(C)$ while the set $X \in H_V$ remains the same:

**Theorem 1c.** *Let $X \in H_V$ to be a subset of a unit hypercube $H$ vertexes. For any $x_k \in X$ the set of hypercube $H(C)$ points $M_k(C) = \{x \in H(C) : x_k\, x \geq Xx\}$ is a convex hull of $H_k(C) = \{h \in H_V(C) : x_k\, h \geq Xh\}$:*

$$M_k(C) = conv(H_k(C))$$

**Proof.** Let $X(C)$ to be a set of $H(C)$ vertexes induced by $X \in H_V$. A hypercube $H(C)$ has a unit volume and integer valued vertexes. Thus, theorem 1b is valid for the set $X(C)$ and a hypecube $H(C)$ and, therefore, $M'_k(C) = conv(H'_k(C))$. From lemma 2.3.1 it follows that $M_k(C) = M'_k(C)$ and $H_k(C) = H'_k(C)$. Therefore, $M_k(C) = M'_k(C) = conv(H'_k(C)) = conv(H_k(C))$. ∎

Theorem 1c provides basis for optimality conditions for a COP in case $c \in H(C)$:

**Theorem 2.1b (extended).** *Let $X \in H \in R^N$ to be a set of incidence vectors induced by a COP and negative components of a weight vector $c$ are indexed by $C \in E$. Let $H_k(C) \in H_V(C)$ to be a subset of $H(C)$ vertexes defined by $x_k \in X$ so that $H_k(C) = \{h \in H_V(C) : x_k\, h \geq Xh\}$. Then $x_k \in X$ is optimal if and only if $c \in cone(H_k(C))$:*

$$\exists\, \gamma_j \geq 0 \text{ so that } c = \sum\nolimits_{h_j \in H_k(C)} \gamma_j h_j \qquad (2)$$

**Proof.** See theorem 1c. ∎

Theorem 2.1b states that the set $M_k(C) = \{c \in H(C) : x_k\, c \geq Xc\}$ is defined by a convex hull of the hypercube $H(C)$ vertexes: $H_k(C) = \{h \in H_V(C) : x_k\, h \geq Xh\}$. Let us extend theorem 2.1b for a general case $c \in H^{\pm 1}$ (or $c \in R^N$ if we do not scale it):

**Theorem 2.1c (general).** *Let $X \in H \in R^N$ to be a set of incidence vectors induced by a COP and $c \in R^N$. Let $V$ to be a set of all -1/0/1 valued vectors in $R^N$. Let $H_k(V) \in V$ to be a subset of $V$ defined by $x_k \in X$ so that $H_k(V) = \{h \in V : x_k\, h \geq Xh\}$. Then $x_k \in X$ is optimal if and only if $c \in cone(H_k(V))$:*

$$\exists\, \gamma_j \geq 0 \text{ so that } c = \sum\nolimits_{h_j \in H_k(V)} \gamma_j h_j \qquad (3)$$

**Proof.** Let us define the set $M_k(H^{\pm 1}) = \{c \in H^{\pm 1} : x_k c \geq Xc\}$ of weigth vector $c \in H^{\pm 1}$ values so that $x_k \in X$ is optimal. As $H^{\pm 1} = \cup_{C \in E} H(C)$ it is obvious that $M_k(H^{\pm 1}) = \cup_{C \in E} M_k(C)$ and, by theorem 2.1b, we have $\cup_{C \in E} M_k(C) = \cup_{C \in E} conv(H_k(C)) = conv(\cup_{C \in E} H_k(C))$ where $H_k(C) = \{h \in H_V(C) : x_k h \geq Xh\}$. The set $\cup_{C \in E} H_V(C)$ is just a set of all -1/0/1 valued vectors in $R^N$ and, thus, $\cup_{C \in E} H_k(C) = H_k(V) = \{h \in V : x_k h \geq Xh\}$. Which, in turn, proves that $M_k(H^{\pm 1}) = \cup_{C \in E} M_k(C) = \cup_{C \in E} conv(H_k(C)) = conv(\cup_{C \in E} H_k(C)) = conv(H_k(V))$. So in case $c \in R^N$ the incident vector $x_k \in X$ is optimal if and only if $c \in cone(H_k(V))$. ∎

## 3. Complete description of polytopes in combinatorial optimization

A general powerful method of proving the validity of an inequality was proposed by (Chvatal 1973) and is essentially the same as the cut generation method of Gomory (Gomory 1963). It is based upon the fact that there is often an integer programming formulation of combinatorial optimization problem. All valid inequalities for the convex hull of incidence vectors can be derived from this (Pulleyblank 1989). Sometimes, it is easier to derive a linear description of some high dimensional polytope by introducing extra variables to formulation of original COP (i.e. variables that do not correspond to groundset $E$). Such formulation is commonly referred to as extended. In some cases extended formulation can be used to find description of the polytope for original COP (Schrijver 2004).

Theorem 2.1c provides a basis for a different procedure to get a complete set of valid and non-redundant facet inducing inequalities for any COP even without analytical integer programing formulation or extended formulation. Any COP is transformed into linear program over the polytope $P^X$, however, in many cases without a complete description of $P^X$ as a system of linear inequalities. From complementary slackness theorem it follows that optimality conditions for any linear program at an extreme point $x_k$ is defined by a complete set of facets that intersects at $x_k$. We proved that optimality conditions (3) is necessary and sufficient, therefore, the set of independent cone generators $\widetilde{H}_k(V) = \{h \in H_k(V) : h \notin cone(H_k(V)\setminus\{h\})\}$ for $cone(H_k(V))$ represent a complete nonredundant set of $P^X$ facets that intersects at $x_k$. This provide a basis for complete description of the polytopes induced by combinatorial optimization problems:

**Theorem 3 (complete description of polytopes in combinatorial optimization).** *Let $X \in R^N$ to be a set of incidence vectors induced by a COP. Let $P^X = conv(X)$ to be the induced polytope. Let $V \in R^N$ to be a set of all -1/0/1 valued vectors. Let $H_k(V)$ to be a subset of $V$ defined for each $x_k \in X$ so that $H_k(V) = \{h \in V : x_k h \geq Xh\}$ and $\widetilde{H}_k(V) = \{h \in H_k(V) : h \notin cone(H_k(V)\setminus\{h\})\}$ to be a set of independent cone generators for $cone(H_k(V))$. Then $P_V^X = \cup_{x_k \in X} \widetilde{H}_k(V)$ defines the complete set of nonredundant facet inducing inequalities for $P^X$:*

*inequality $hx \leq l$ induces a nonredundant facet for each $h \in P_V^X$ and $l = hx_k$ for $x_k \in X$ so that $h \in \widetilde{H}_k(V)$*

**Proof.** The proof of theorem 3 directly follows from complementary slackness theorem and theorem 2.1c. ∎

Let us remark that for symmetric problems such as traveling salesman problem (TSP) it is enough to deduce facets only for one extreme point as all facets for the other extreme points can be deduced by proper permutation of corresponding indexes. For combinatorial problems that have the same sizes of all feasible subsets ( $|t| = n$ for all $t \in T$) it is sufficient to consider only case when $c \in H$ as we can scale $c$ and add/subtract unit vector to/from $c$ without affecting hierarchical order of solutions induced by $c$ (TSP is also from this class). And even if some components of $c$ are negative by proper transformation (scaling + adding unit vector) we can ensure that $c \in H$. Therefore, for this class of problems the set $P_V^X$ from theorem 3 can be redefined as $P_H^X$ using only a unit hypercube vertexes $H_V$ instead of $V$ (except non-negativity constraints):

$$P_H^X = \{h \in H_V : h \in H_k \text{ and } h \notin cone(H_k \setminus \{h\}) \text{ for at least one } x_k \in X\}$$

## 4. Final remarks

In this paper we have provided analytical formulation of optimality conditions for a general linear combinatorial optimization problem for the first time up to our knowledge. This provides basis for a complete description of any polytope induced by combinatorial optimization (as a system of linear inequalities) even without integer programming formulation (or extended formulation).

## References

Aardal, K. and S. v. Hoesel (1996). "Polyhedral techniques in combinatorial optimization." Statistica Neerlandica **50**: 3-26.

Chvatal, V. (1973). "Edmonds polytopes and a hierarchy of combinatorial problems " Discrete Mathematics **4**: 305-337.

Gomory, R. E. (1963). "An algorithm to integer solutions to linear programs." Resent advances in Mathematical Programing: 269-302.

Pulleyblank, W. R. (1989). Polyhedral combinatorics, Optimization (G.L. Nemhauser, A.H.G. Rinnooy Kan, and M.J. Todd, eds.). Handbooks in OR & MS. **1:** 371-446.

Schrijver, A. (2004). Combinatorial optimization (Polyhedra ans Efficiency) Springer. **A-C**.